# MATROID BASE POLYTOPE DECOMPOSITION

V. CHATELAIN AND J. L. RAMÍREZ ALFONSÍN

ABSTRACT. Let $P(M)$ be the matroid base polytope of a matroid $M$. A *matroid base polytope decomposition* of $P(M)$ is a decomposition of the form $P(M) = \bigcup_{i=1}^{t} P(M_i)$ where each $P(M_i)$ is also a matroid base polytope for some matroid $M_i$, and for each $1 \leq i \neq j \leq t$, the intersection $P(M_i) \cap P(M_j)$ is a face of both $P(M_i)$ and $P(M_j)$. In this paper, we investigate *hyperplane splits*, that is, polytope decompositions when $t = 2$. We give sufficient conditions for $M$ so $P(M)$ has a hyperplane split and characterize when $P(M_1 \oplus M_2)$ has a hyperplane split where $M_1 \oplus M_2$ denote the *direct sum* of matroids $M_1$ and $M_2$. We also prove that $P(M)$ has not a hyperplane split if $M$ is binary. Finally, we show that $P(M)$ has not a decomposition if its 1-skeleton is the *hypercube*.

## 1. INTRODUCTION

For general background in matroid theory we refer the reader to [18, 21]. A *matroid* $M = (E, \mathcal{B})$ of *rank* $r$ is a finite set $E = \{1, \ldots, n\}$ (called the *ground* set of $M$) together with a nonempty collection $\mathcal{B} = \mathcal{B}(M)$ of $r$-subsets of $E$ (called the *bases* of $M$) satisfying the following *basis exchange axiom*:

if $B_1, B_2 \in \mathcal{B}$ and $e \in B_1 \setminus B_2$ then there exists $f \in B_2 \setminus B_1$ such that $(B_1 - e) + f \in \mathcal{B}$.

The *independent* set of $M$, denoted by $\mathcal{I}(M)$ is given by all subsets of bases of $M$. For a matroid $M = (E, \mathcal{B})$, let $P(M)$ be the *matroid base polytope* of $M$ defined as the convex hull of the incidence vectors of bases of $M$, that is,

$$P(M) := \text{conv} \left\{ \sum_{i \in B} e_i : B \text{ a base of } M \right\},$$

where $e_i$ denotes the $i^{th}$ standard basis vector in $\mathbb{R}^n$. $P(M)$ is a polytope of dimension at most $n - 1$. Notice that $P(M)$ is a face of the *independent set polytope* $I(M)$ which is obtained as the convex hull of the incidence vectors of the independent sets of $M$. These polytopes were first studied by Edmonds [5, 6].

A *matroid base polytope decomposition* of $P(M)$ is a decomposition

$$P(M) = \bigcup_{\substack{i=1 \\ 1}}^{t} P(M_i)$$

where each $P(M_i)$ is also a matroid base polytope for some matroid $M_i$, and for each $1 \leq i \neq j \leq t$, the intersection $P(M_i) \cap P(M_j)$ is a face of both $P(M_i)$ and $P(M_j)$.

$P(M)$ is said to be *decomposable* if it has a matroid base polytope decomposition with $t \geq 2$, and *indecomposable* otherwise. A decomposition is called *hyperplane split* if $t = 2$.

Matroid base polytope decomposition have appeared in many different contexts. For instance, they are treated in the work by Hacking, Keel and Tevelev [9, Section 3.3] in relation with the compactification of the moduli space of hyperplane arrangements (see also [10] and [11, Section 2.6]), in Speyer's work [19, 20] concerning tropical linear spaces, and in Lafforgue's work [12, 13] while studying the compactifications of the fine Schubert cell of the Grassmannian. In particular, Lafforgue's work implies that for a matroid $M$ represented by vectors in $\mathbb{F}^r$, if $Q(M)$ is indecomposable, then $M$ will be *rigid*, that is, $M$ will have only finitely many realizations, up to scaling and the action of $GL(r, \mathbb{F})$. Recently, Billera, Jia and Reiner [3] (see also the closely related results due to Luoto [14]), Speyer [19, 20], Derksen [4] and Ardila, Fink and Rincon [1] have showed that different interesting matroids functions (as quasisymmetric functions and Tutte's polynomials) behaves like *valuations* on the the associated matroid base polytope decomposition.

It is therefore of interest to know whether a given matroid base polytope is decomposable or not. Unfortunately, there is not much known about the existence or nonexistence of such decompositions (even for the case $t = 2$). Kapranov [10, Section 1.3] showed that all decompositions of a (convenient parametrizied) rank two matroid can be achieved by a sequence of hyperplane splits. In [3], Billera, Jia and Reiner found five rank three matroids on 6 elements for which the corresponding polytope is indecomposable. They also showed that $P(M)$ can be splited into three indecomposable pieces where $M$ is the rank three matroid on $\{1, \ldots, 6\}$ having every triple but $\{1, 2, 3\}, \{1, 4, 5\}$ and $\{3, 5, 6\}$ as bases. Moreover, they showed that this decomposition cannot be obtained via hyperplane splits. In this paper, we show the existence and nonexistence of hyperplane splits for some infinite classes of matroids. We also give a special family of matroid base polytopes that are indecomposable.

It is known that nonempty faces of matroid base polytope are matroid base polytopes [8, Theorem 2]. So, the common face $P(M_i) \cap P(M_j)$ (whose vertices correspond to elements of $\mathcal{B}(M_i) \cap \mathcal{B}(M_j)$) must also be a matroid base polytope. Thus, in order to investigate the nonexistence of base polytope decomposition, one may consider the following combinatorial decomposition version. A *matroid base decomposition* of a matroid $M$ is a decomposition

$$\mathcal{B}(M) = \bigcup_{i=1}^{t} \mathcal{B}(M_i)$$

where $\mathcal{B}(M_k)$, $1 \leq k \leq t$ and $\mathcal{B}(M_i) \cap \mathcal{B}(M_j)$, $1 \leq i \neq j \leq t$ are collections of bases of matroids.

$M$ is said to be *combinatorial decomposable* if it has a matroid base decomposition. We say that the decomposition is nontrivial if $\mathcal{B}(M_i) \neq \mathcal{B}(M)$ for all $i$. If $P(M)$ is decomposable then $M$ is clearly combinatorial decomposable. However, a matroid base decomposition do not necessarily induce a matroid base polytope decomposition. For instance, the rank 2 matroid $M^*$ where $\mathcal{B}(M^*) = \{\{1,2\},\{1,3\},\{2,3\},\{2,4\},\{3,4\}\}$ has a combinatorial decomposition given by $\mathcal{B}(M_1) = \{\{1,2\},\{2,3\},\{2,4\}\}$ and $\mathcal{B}(M_2) = \{\{1,3\},\{2,3\},\{3,4\}\}$ since $\mathcal{B}(M_1), \mathcal{B}(M_2)$ and $\mathcal{B}(M_1) \cap \mathcal{B}(M_2) = \{2,3\}$ are collections of bases of matroids. However the corresponding polytopes $P(M_1)$ and $P(M_2)$ do not decompose $P(M^*)$, see Figure 1 (a).

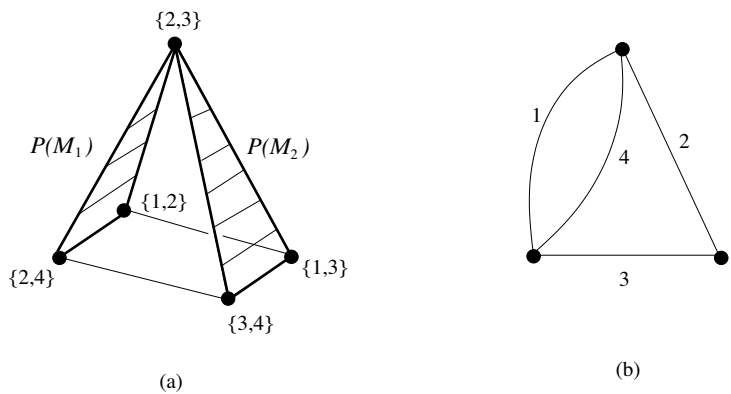

FIGURE 1. (a) $P(M^*), P(M_1)$ and $P(M_2)$ (b) A graph $G$ where each spanning trees correspond to a base of $M^*$.

We notice that $M^*$ is actually a *graphic* matroid (a matroid $M(E, \mathcal{B})$ is said to be *graphic* if there exists a graph $G$ with $|E|$ edges such that there exists a bijection between the elements $E$ of $M$ and the edges of $G$ so that the elements of each base of $M$ correspond to a set of edges of a spanning forest of $G$). The required graph $G$ for $M^*$ is given in Figure 1 (b).

In order to prove the existence of nontrivial hyperplane splits, we always first show the existence of matroid base decompositions and then prove that this induces a hyperplane split since some geometric conditions (see Proposition 1)are verified.

In next section, we shall give sufficient conditions for a matroid $M$ so that $P(M)$ has a nontrivial hyperplane split. Our constructive method allows to show the existence of at least $\lfloor \frac{n}{2} \rfloor$ *different* hyperplane splits of $P(U_{n,r})$ with $n \geq r + 2 \geq 3$ where $U_{n,r}$ denotes the *uniform matroid* on $n$ elements of rank $r$ (recall that $\mathcal{B}(U_{n,r})$ consist of all $r$-subsets of $\{1, \ldots, n\}$). In Section 3, we present a complete characterization for matroid $M_1 \oplus M_2$ so $P(M_1 \oplus M_2)$ has a nontrivial hyperplane split where $M_1 \oplus M_2$ denote the *direct sum* of matroids $M_1$ and $M_2$. In Section 4, we will show that $P(M)$ has not a nontrivial hyperplane split if $M$ is a *binary* matroid, that is, if $M$ is *representable* over $\mathbb{F}^2$. We finally prove that if the 1-skeleton of $P(M)$ is the hypercube then $P(M)$ is indecomposable.

## 2. Decomposition

The *matroid base graph* $G(M)$ of a matroid $M$ is the graph having as set of vertices the bases of $M$ and there is an edge between two vertices (bases) $B_1, B_2$ if and only if there exist a pair of elements $e \in B_1$ and $f \in B_2$ such that $B_2 = (B_1 \setminus e) + f$, that is, the symmetric difference of $B_1$ and $B_2$, denoted by $\Delta(B_1, B_2)$, is equals two. It is known [7] that $G(M)$ is the 1-skeleton of $P(M)$ (in other words, the edges of $P(M)$ represent the basis exchange axiom) and that $G(M)$ is connected. We present the following geometric result used throughout the rest of the paper.

**Proposition 1.** *Let $P$ be a $d$-polytope with set of vertices $X$. Let $H$ be a hyperplane such that $H \cap P \neq \emptyset$ with $H$ not supporting $P$. So, $H$ splits $P$ into polytopes $P_1$ and $P_2$, that is, $H \cap P = P_1 \cap P_2 = F \neq \emptyset$. $H$ also partition $X$ into sets $X_1$ and $X_2$ with $X_1 \cap X_2 = W$. Then, for each edge $[u, v]$ of $P$ we have that $\{u, v\} \subset X_i$ with either $i = 1$ or 2 if and only if $F = conv(W)$.*

*Proof.* We notice that $X_1, X_2 \neq \emptyset$ (since $H$ is not supporting). and let $[u, v]$ be an edge of $P$.

(*Necessity*) We shall proceed by contradiction. Suppose that $\{u, v\} \subset X_i$ with $i = 1$ or 2 and that $F \neq conv(W)$. Since $conv(W) \subset F$ then there exists a vertex $x$ in $F$ such that $x \notin W$ (and thus $x \notin X$). So, $x$ is the intersection of $H$ with an edge $[u, v]$ of $P$ with $u \in X_1 \setminus X_2$ and $v \in X_2 \setminus X_1$ which is a contradiction.

(*Sufficiency*) We shall proceed by contradiction. Suppose that $F = conv(W)$ and that $u \in X_1 \setminus X_2$ and $v \in X_2 \setminus X_1$. So, $u \in P_1 \setminus P_2$ and $v \in P_2 \setminus P_1$ and therefore $[u, v] \cap F = s$ with $s$ a vertex of $F$ different of $u$ and $v$. Then, $conv(W)$ does not contain $s$ (since $s$ is not a vertex of $W$) and so $F \neq conv(W)$, which is a contradiction.    □

We have the following easy consequence of this proposition.

**Corollary 1.** *Let $P, P_1$ and $P_2$ be the polytopes as in Proposition 1. Then, $F = conv(W)$ if and only if $P_i = conv(X_i)$, $i = 1, 2$ (and thus $P = P_1 \cup P_2$ with $P_1$ and $P_2$ polytopes of the same dimension as $P$ sharing a facet).*

Let $M = (E, \mathcal{B})$ be a matroid of rank $r$ and let $A \subseteq E$. We recall that the independent set of the *restriction* matroid of $M$ to $A$, denoted by $M|_A$, is given by $\mathcal{I}(M|_A) = \{I \subseteq A : I \in \mathcal{I}(M)\}$.

Let $(E_1, E_2)$ be a partition of $E$, that is, $E = E_1 \cup E_2$ and $E_1 \cap E_2 = \emptyset$. Let $r_i > 1$, $i = 1, 2$ be the rank of $M|_{E_i}$. We say that $(E_1, E_2)$ is a *good* partition if there exist integers $0 < a_1 < r_1$ and $0 < a_2 < r_2$ with the following properties:

(P1) $r_1 + r_2 = r + a_1 + a_2$ and

(P2) for all $X \in \mathcal{I}(M|_{E_1})$ with $|X| \leq r_1 - a_1$ and all $Y \in \mathcal{I}(M|_{E_2})$ with $|Y| \leq r_2 - a_2$ we have $X \cup Y \in \mathcal{I}(M)$.

**Lemma 1.** *Let $M = (E, \mathcal{B})$ be a matroid of rank $r$ and let $(E_1, E_2)$ be a good partition of $E$. Let*

$$\mathcal{B}(M_1) = \{B \in \mathcal{B}(M) : |B \cap E_1| \leq r_1 - a_1\}$$

*and*

$$\mathcal{B}(M_2) = \{B \in \mathcal{B}(M) : |B \cap E_2| \leq r_2 - a_2\}.$$

*where $r_i$ is the rank of matroid $M|_{E_i}$, $i = 1, 2$ and $a_1, a_2$ are integers verifying Properties (P1) and (P2). Then, $\mathcal{B}(M_1)$ and $\mathcal{B}(M_2)$ are the collections of bases of matroids $M_1$ and $M_2$ respectively.*

*Proof.* We shall prove that $\mathcal{B}(M_1)$ is the collection of bases of a matroid (this can be done similarly for $\mathcal{B}(M_2)$). We show that the elements in $\mathcal{B}(M_1)$ verify the basis exchange axiom. Let $X, Y \in \mathcal{B}(M_1) \subset \mathcal{B}(M)$ and suppose that $e \in X \setminus Y$. Since $M$ is a matroid then there exists $f \in Y \setminus X$ such that $X - e + f \in \mathcal{B}(M)$. We have two cases.

Case 1) Suppose that either $|X \cap E_1| < r_1 - a_1$ or $e \in E_1$. Then, $|(X - e + f) \cap E_1| \leq r_1 - a_1$ and thus $X - e + f \in \mathcal{B}(M_1)$.

Case 2) Suppose that $|X \cap E_1| = r_1 - a_1$ and $e \in E_2$. In one hand, we have that $|(X - e) \cap E_2| = r - 1 - (r_1 - a_1) = r - r_1 + a_1 - 1$. On the other hand, we have that $|Y \cap E_1| \leq r_1 - a_1$ and thus $|Y \cap E_2| \geq r - (r_1 - a_1)$. So, there exists $g \in Y \setminus (X - e) = Y \setminus X$ with $g \in E_2$ such that $(X - e + g) \cap E_2 \in \mathcal{I}(M|_{E_2})$ and also, since $g \notin E_1$, then $(X - e + g) \cap E_1 \in \mathcal{I}(M|_{E_1})$. Moreover, $|(X - e + g) \cap E_1| = r_1 - a_1$ and $|(X - e + g) \cap E_2| = r - r_1 + a - 1 = r_2 - a_2$ and thus, by Property (P1), $|X - e + g| = r_1 + r_2 - a_1 - a_2 = r$. So, by Property (P2), we have $((X - e + g) \cap E_1) \cup (X - e + g) \cap E_2)) = X - e + g \in \mathcal{I}(M)$. But since, $|X - e + g| = r$ then $X - e + g \in \mathcal{B}(M)$ and since $|(X - e + g) \cap E_1| = r_1 - a_1$ then $X - e + g \in \mathcal{B}(M_1)$. $\quad\square$

Notice that property (P2), needed for the proof of Lemma 1, can be replaced by the following weaker condition.

(P2') for all $B_1 \in \mathcal{B}(M_1)$ and all $Y \in \mathcal{I}(M|_{E_2})$ with $|Y| \leq r_2 - a_2$ we have $(B_1 \cap E_1) \cup Y \in \mathcal{I}(M)$ and for all $B_2 \in \mathcal{B}(M_2)$ and all $Y \in \mathcal{I}(M|_{E_1})$ with $|Y| \leq r_1 - a_1$ we have $(B_2 \cap E_2) \cup Y \in \mathcal{I}(M)$.

(P2) is clearly a stronger condition that (P2') since

$$\{B_1 \cap E_1 : B_1 \in \mathcal{B}(M_1)\} \subset \{X : |X| \leq r_1 - a_1, X \in \mathcal{I}(M|_{E_1})\}$$

and

$$\{B_2 \cap E_2 : B_2 \in \mathcal{B}(M_2)\} \subset \{X : |X| \leq r_2 - a_2, X \in \mathcal{I}(M|_{E_2})\}.$$

**Theorem 1.** *Let $M = (E, \mathcal{B})$ be a matroid of rank $r$ and let $(E_1, E_2)$ be a good partition of $E$. Let $M_1$ and $M_2$ be matroids given in Lemma 1. Then, $P(M) = P(M_1) \cup P(M_2)$ is a nontrivial hyperplane split.*

*Proof.* We will first show that $\mathcal{B}(M_1)$ and $\mathcal{B}(M_2)$ give a nontrivial matroid base decomposition of $M$. For this, we show

$(i)$ $\mathcal{B}(M) = \mathcal{B}(M_1) \cup \mathcal{B}(M_2)$,

$(ii)$ $\mathcal{B}(M_1), \mathcal{B}(M_2) \subset \mathcal{B}(M)$,

$(iii)$ $\mathcal{B}(M_1), \mathcal{B}(M_2) \nsubseteq \mathcal{B}(M_1) \cap \mathcal{B}(M_2)$,

$(iv)$ $\mathcal{B}(M_1) \cap \mathcal{B}(M_2) \neq \emptyset$ and

$(v)$ $\mathcal{B}(M_1) \cap \mathcal{B}(M_2)$ is the collection of bases of a matroid.

We then show that this matroid base decomposition induces a nontrivial hyperplane split. For this, we show

$(vi)$ there exists an hyperplane containing the vertices corresponding to $\mathcal{B}(M_1) \cap \mathcal{B}(M_2)$, and not supporting $P(M)$,

$(vii)$ any edge of $P(M)$ is also and edge of either $P(M_1)$ or $P(M_2)$.

So, by Corollary 1, $\mathcal{B}(M_1) \cap \mathcal{B}(M_2)$ are the set of vertices of a facet of $P(M_1)$ and $P(M_2)$ and thus $P(M) = P(M_1) \cup P(M_2)$.

We may now prove the above claims.

$(i)$ We claim that $\overline{\mathcal{B}(M)} = \emptyset$. Indeed, let $B \in \overline{\mathcal{B}(M)} = \{B \in \mathcal{B}(M) : |B \cap E_1| > r_1 - a_1$ and $|B \cap E_2| > r_2 - a_2\}$. Since $E_1 \cap E_2 \neq \emptyset$ then $|B| > r_1 + r_2 - a_1 - a_2 = r$ which is not possible. So, $\mathcal{B}(M) = \mathcal{B}(M_1) \cup \mathcal{B}(M_2)$.

$(ii)$ We show that $\mathcal{B}(M_1) \subset \mathcal{B}(M)$ (it can also be proved that $\mathcal{B}(M_2) \subset \mathcal{B}(M)$ by using similar arguments). Let $B$ be a base of $M|_{E_1}$ (and so $|B| = r_1$). We have two cases.

Case 1) If $r_1 = r$ then $B \in \mathcal{B}(M)$ but $B \notin \mathcal{B}(M_1)$ since $|B \cap E_1| = r > r_1 - a_1$.

Case 2) If $r_1 < r$ then there exists a base $B' \in \mathcal{B}(M)$ with $B \subset B'$ since $B \in \mathcal{I}(M)$. Moreover, $|B' \cap E_1| \geq r_1 > r_1 - a_1$ and thus $B' \in \mathcal{B}(M)$ but $B' \notin \mathcal{B}(M_1)$.

$(iii)$ If $\mathcal{B}(M_1) \subseteq \mathcal{B}(M_1) \cap \mathcal{B}(M_2)$ then $\mathcal{B}(M_1) \subseteq \mathcal{B}(M_2)$ and thus $\mathcal{B}(M_2) = \mathcal{B}(M_1) \cup \mathcal{B}(M_2) = \mathcal{B}(M)$ contradicting $(ii)$.

$(iv)$ Let $X, Y \in \mathcal{B}(M)$ such that $|X \cap E_1| = r_1$ and $|Y \cap E_2| = r_2$ (we have seen in $(ii)$ that such bases always exist). Since $G(M)$ is connected then there exists a path $X = B_1, \ldots, B_m = Y$ connecting $X$ and $Y$ where $B_i \in \mathcal{B}(M)$ for each $i = 1, \ldots, m$. Since $|B_i \Delta B_{i+1}| = 2$ then $|B_i \cap E_1| = |B_{i+1} \cap E_1| + k$ with either $k = -1, 0, +1,$. Moreover, $|B_1 \cap E_1| = r_1$ and $|B_m \cap E_1| = r - r_2 = r_1 - a_1 - a_2$. So, there must exists an index $k$ such that $|B_k \cap E_1| = r_1 - a_1$. Moreover, $|B_k \cap E_2| = r - r_1 + a_1 = r_2 - a_2$. So, $B_k \in \mathcal{B}(M_1) \cap \mathcal{B}(M_2)$.

$(v)$ This can be done by using similar arguments as those used in Lemma 1.

$(vi)$ We first notice that if $B \in \mathcal{B}(M_1) \cap \mathcal{B}(M_2)$ then $B \in \mathcal{B}(M)$ (so $|B| = r$) and $|B \cap E_1| \leq r_1 - a_1$, $|B \cap E_2| \leq r_2 - a_2 = r - r_1 + a_1$. So, if $B \in \mathcal{B}(M_1) \cap \mathcal{B}(M_2)$ then $|B \cap E_1| = r_1 - a_1$ and $|B \cap E_2| = r_2 - a_2$. Let $H$ be the hyperplane defined by

$$H(\mathbf{x}) = \left\{ \mathbf{x} \in \mathbb{R}^n : \langle \mathbf{e}, \mathbf{x} \rangle = r_1 - a_1 \text{ where } \mathbf{e} = \sum_{i \in E_1} e_i \right\}.$$

If $B \in \mathcal{B}(M_1)$ then $|B \cap E_1| \leq r_1 - a_1$ and thus $H(x_B) \leq r_1 - a_1$. If $B \in \mathcal{B}(M_2)$ then $|B \cap E_2| \leq r_2 - a_2$ and so $|B \cap E_1| \geq r - r_2 + a_2 = r_1 - a_1$. Then, $H(x_B) \geq r_1 - a_1$. Therefore, $H$ contains only elements in $\mathcal{B}(M_1) \cap \mathcal{B}(M_2)$. Notice that $dim(P(M) \cap H) = d - 1$ where $d = dim(P(M))$. Indeed, $dim(P(M) \cap H) \leq d$ and if $dim(P(M) \cap H) \leq d - 2$ then $P(M)$ would be contained in the closed-half space $\overline{H}^+$ or in $\overline{H}^-$ which is impossible by $(iii)$. So, the affine space $H'$ containing $P(M) \cap H$ is a hyperplane in $\mathbb{R}^d$ not supporting $P(M)$.

$(vii)$ Suppose that there exists and edge in $P(M)$ belonging to neither $P(M_1)$ nor $P(M_2)$. So, there exist $B_1 \in \mathcal{B}(M_1) \backslash \mathcal{B}(M_2)$ and $B_2 \in \mathcal{B}(M_2) \backslash \mathcal{B}(M_1)$ with $\Delta(B_1, B_2) = 2$. Then, $|B_1 \cap E_1| < r_1 - a_1$ and $|B_2 \cap E_1| > r_1 - a_1$. Combining the latter with the fact that $\Delta(B_1, B_2) = 2$, then we must have $|B_1 \cap E_1| = r_1 - a_1$ which is a contradiction. □

**Remark 1.** *It is not necessarily to have a good partition in the hypothesis of Theorem 1 but only the existence of a partition $(E_1, E_2)$ of $E$ verifying property (P1) and such that the sets $\mathcal{B}(M_1)$, $\mathcal{B}(M_2)$ and $\mathcal{B}(M_1) \cap \mathcal{B}(M_2)$ be collections of bases of matroids. Indeed, this is the case if the partition verifies either (P2) or (P2').*

**Example 1:** Let us consider $U_{4,2}$ and the good partition $E_1 = \{1, 2\}$ and $E_2 = \{3, 4\}$ (and so $r_1 = r_2 = 2$) with $a_1 = a_2 = 1$. Then $\mathcal{B}(M_1) = \{\{1, 3\}, \{1, 4\}, \{2, 3\}, \{2, 4\}, \{3, 4\}\}$, $\mathcal{B}(M_2) = \{\{1, 2\}, \{1, 3\}, \{1, 4\}, \{2, 3\}, \{2, 4\}\}$ and $\mathcal{B}(M_1) \cap \mathcal{B}(M_2) = \{\{1, 3\}, \{1, 4\}, \{2, 3\}, \{2, 4\}\}$, see Figure 2.

We remark that there exist matroids with ground set $E$ not having a good partition. For instance, let $M(K_4)$ be the matroid associated to the complete graph on four vertices. So, the bases of $M(K_4)$ are all triples except $\{1, 2, 5\}$, $\{2, 3, 6\}$, $\{3, 4, 5\}$ and $\{1, 4, 6\}$ (corresponding to circuits of $K_4$), see Figure 3. Since $1 < r_1, r_2 \leq 3$ then $2 \leq |E_1|, |E_2| \leq 4$. We only have the following four possibilities (up to symmetries), these are illustrated in Figure 3 (a)-(d).

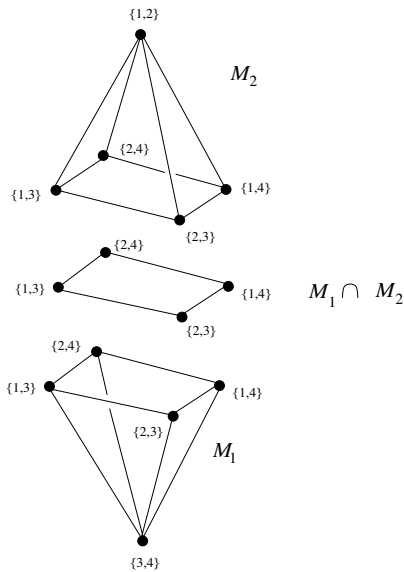

FIGURE 2. Hyperplane split of $P(U_{4,2})$

(a) If $E_1 = \{1,2,3\}$ and $E_2 = \{4,5,6\}$ then $r_1 = 3, r_2 = 3$ and $\mathcal{B}(M_1) = \{B \in \mathcal{B}(M) : |B \cap E_1| \leq 2\}$, $\mathcal{B}(M_2) = \{B \in \mathcal{B}(M) : |B \cap E_2| \leq 1\}$. We have that $\mathcal{B}(M_2)$ is not the collection of bases of a matroid. Indeed, $\{1,3,5\}, \{1,2,4\} \in \mathcal{B}(M_2)$ so, by the basis exchange axiom, we have that $\{1,2,5\}$ or $\{1,4,5\}$ should be in $\mathcal{B}(M_2)$. But, $\{1,2,5\} \notin \mathcal{B}(M(K_4))$ (and so $\{1,2,5\} \notin \mathcal{B}(M_2)$) and $\{1,4,5\} \cap E_2 = \{4,5\}$ (and so $\{1,4,5\} \notin \mathcal{B}(M_2)$).

(b) If $E_1 = \{1,2,5\}$ and $E_2 = \{3,4,6\}$ then $r_1 = 2, r_2 = 3$ and $\mathcal{B}(M_1) = \{B \in \mathcal{B}(M) : |B \cap E_1| \leq 1\}$, $\mathcal{B}(M_2) = \{B \in \mathcal{B}(M) : |B \cap E_2| \leq 2\}$. We have that $\mathcal{B}(M_1)$ is not the collection of bases of a matroid. Indeed, $\{1,3,4\}, \{2,4,6\} \in \mathcal{B}(M_1)$ so, by the basis exchange axiom, we have that $\{1,2,4\}$ or $\{1,4,6\}$ should be in $\mathcal{B}(M_1)$. But, $\{1,4,6\} \notin \mathcal{B}(M(K_4))$ (and so $\{1,4,6\} \notin \mathcal{B}(M_1)$) and $\{1,2,4\} \cap E_1 = \{1,2\}$ (and so $\{1,2,4\} \notin \mathcal{B}(M_1)$).

(c) If $E_1 = \{1,2\}$ and $E_2 = \{3,4,5,6\}$ then $r_1 = 2, r_2 = 3$ and $\mathcal{B}(M_1) = \{B \in \mathcal{B}(M) : |B \cap E_1| \leq 1\}$, $\mathcal{B}(M_2) = \{B \in \mathcal{B}(M) : |B \cap E_2| \leq 2\}$. We have that $\mathcal{B}(M_2)$ is not the collection of bases of a matroid. Indeed, $\{2,4,6\}, \{1,3,6\} \in \mathcal{B}(M_2)$ so, by the basis exchange axiom, we have that $\{1,4,6\}$ or $\{3,4,6\}$ should be in $\mathcal{B}(M_2)$. But, $\{1,4,6\} \notin \mathcal{B}(M(K_4))$ (and so $\{1,4,6\} \notin \mathcal{B}(M_2)$) and $\{3,4,6\} \cap E_2 = \{3,4,6\}$ (and so $\{3,4,6\} \notin \mathcal{B}(M_2)$).

(d) If $E_1 = \{1,3\}$ and $E_2 = \{2,4,5,6\}$ then $r_1 = 2, r_2 = 3$ and $\mathcal{B}(M_1) = \{B \in \mathcal{B}(M) : |B \cap E_1| \leq 1\}$, $\mathcal{B}(M_2) = \{B \in \mathcal{B}(M) : |B \cap E_2| \leq 2\}$. We have that $\mathcal{B}(M_2)$

is not the collection of bases of a matroid. Indeed, $\{1,4,5\}, \{3,4,6\} \in \mathcal{B}(M_2)$ so, by the basis exchange axiom, we have that $\{3,4,5\}$ or $\{4,5,6\}$ should be in $\mathcal{B}(M_2)$. But, $\{3,4,5\} \notin \mathcal{B}(M(K_4))$ (and so $\{3,4,5\} \notin \mathcal{B}(M_2)$) and $\{4,5,6\} \cap E_2 = \{4,5,6\}$ (and so $\{4,5,6\} \notin \mathcal{B}(M_2)$).

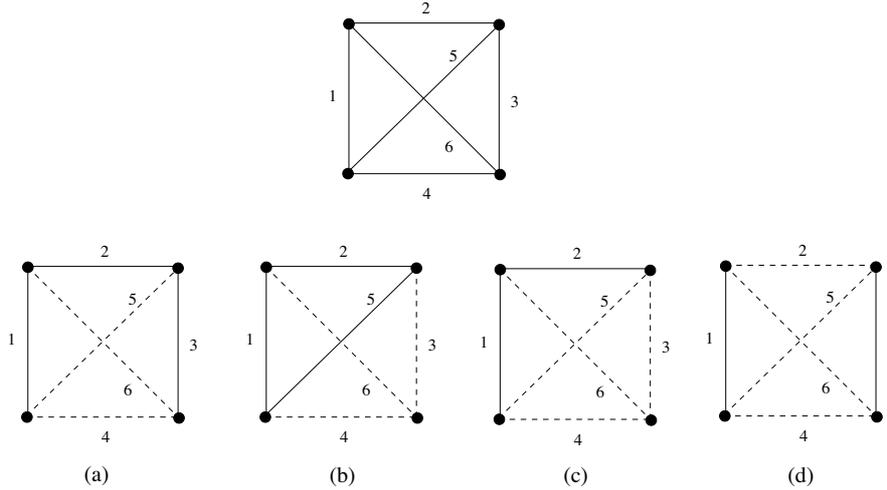

FIGURE 3. $K_4$ and four possibles good partitions. In each case the elements in $E_1$ (resp. in $E_2$) are represented by solid (resp. by dashed) edges.

It turns out that $P(M(K_4))$ has not a nontrivial hyperplane split (this is justified at the end of Section 4).

Let $M = (E, \mathcal{B})$ be a matroid of rank $r$ and let $X \subset E$ be both a circuit and a hyperplane of $M$ (recall that a *hyperplane* is a *flat*, that is $X = cl(X)$, of rank $r - 1$). It is known [18, Proposition 1.5.13] that $\mathcal{B}(M') = \mathcal{B}(M) \cup X$ is the collection of bases of a matroid $M'$ (called, *relaxation* of $M$).

**Corollary 2.** *Let $M = (E, \mathcal{B})$ be a matroid and let $(E_1, E_2)$ be a good partition of $E$. Then, $P(M')$ has a nontrivial hyperplane split where $M'$ is a relaxation of $M$.*

*Proof.* It can be checked that the desired hyperplane split of $P(M')$ can be obtained by using the same given good partition $(E_1, E_2)$ of $E$. □

Although $P(M(K_4))$ has not a nontrivial hyperplane split, base polytopes of relaxations of $M(K_4)$ may have one.

**Example 2:** Let $W^3$ be the matroid of rank 3 on $E = \{1, \ldots, 6\}$ having as set of bases all 3-subsets of $E$ except the triples $\{1,2,5\}$, $\{1,4,6\}$ and $\{2,3,6\}$, see Figure 4. Notice that $W^3$ is a relaxation of $M(K_4)$ (by relaxing circuit $\{3,4,5\}$) but it is not graphic. It can be checked that $E_1 = \{1, 2, 6\}$ and $E_2 = \{3, 4, 5\}$ (and so $r_1 = r_2 = 3$) with $a_1 = 2$ and $a_2 = 1$ is a good partition. In this case we have

$$\begin{aligned}
\mathcal{B}(M_1) = \quad & \{\{1, 3, 4\}, \{1, 3, 5\}, \{1, 4, 5\}, \{2, 3, 4\}, \{2, 3, 5\}, \{2, 4, 5\}, \\
& \{3, 4, 5\}, \{3, 4, 6\}, \{3, 5, 6\}, \{4, 5, 6\}\},
\end{aligned}$$

$$\begin{aligned}
\mathcal{B}(M_2) = \quad & \{\{1, 2, 3\}, \{1, 2, 4\}, \{1, 2, 6\}, \{1, 3, 4\}, \{1, 3, 5\}, \{1, 3, 6\}, \\
& \{1, 4, 5\}, \{1, 5, 6\}, \{2, 3, 4\}, \{2, 3, 5\}, \{2, 4, 5\}, \{2, 4, 6\}, \\
& \{2, 5, 6\}, \{3, 4, 6\}, \{3, 5, 6\}, \{4, 5, 6\}\},
\end{aligned}$$

and

$$\begin{aligned}
\mathcal{B}(M_1) \cap \mathcal{B}(M_2) = \quad & \{\{1, 3, 4\}, \{1, 3, 5\}, , \{1, 4, 5\}, \{2, 3, 4\}, \{2, 3, 5\}, \\
& \{2, 4, 5\}, \{3, 4, 6\}, \{3, 5, 6\}, \{4, 5, 6\}\}.
\end{aligned}$$

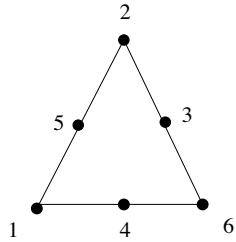

FIGURE 4. Euclidean representation of $W^3$

2.1. **Uniform matroids.** We say that two hyperplane splits $P(M_1) \cup P(M_2)$ and $P(M_1') \cup P(M_2')$ of $P(M)$ are *equivalent* if $P(M_i)$ is *combinatorial equivalent* to $P(M_i')$ for each $i = 1, 2$ (that is, the corresponding face lattices are isomorphic). They are *different* otherwise.

**Corollary 3.** *Let $n \geq r + 2 \geq 4$ be integers. Let $h(U_{n,r})$ be the number of different hyperplane splits of $P(U_{n,r})$. Then,*

$$h(U_{n,r}) \geq \left\lfloor \frac{n}{2} \right\rfloor - 1.$$

*Proof.* For each $k = 2, \ldots, \lfloor \frac{n}{2} \rfloor$, we let $E_1(k) = \{1, \ldots, k\}$ and $E_2(k) = \{k + 1, \ldots, n\}$. So, $M|_{E_1(k)}$ is isomorphic to $U_{k, \min\{k, r\}}$ and $M|_{E_2(k)}$ is isomorphic to $U_{n-k, \min\{n-k, r\}}$. Let $r_1$ and $r_2$ be the ranks of $M|_{E_1(k)}$ and $M|_{E_2(k)}$ respectively. Then, an easy analysis shows that

$$r_1 + r_2 = \min\{n, k + r, n - k + r, 2r\} \geq r + 2.$$

So, we can find integers $a_1, a_2 \geq 1$ such that $r_1 + r_2 = r + a_1 + a_2$ and thus $(E_1(k), E_2(k))$ is a good partition. □

Notice that there might be several choices for the values of $a_1$ and $a_2$ (each of which arises a good partition). However, it is not clear if these partitions give different hyperplane splits.

**Example 3:** Let us consider $U_{n,3}$ with $n \geq 6$. We take the good partition $E_1(k) = \{1, \ldots, k\}$ and $E_2(k) = \{k+1, \ldots, n\}$ for any $2 \leq k \leq \lfloor \frac{n}{2} \rfloor$ (and thus $r_1 = 2, r_2 = 3$). If we set $a_2 = 1$ and $a_1$ such that $r_1 - a_1 = 1$ then $\mathcal{B}(M_1) = \{B \in \mathcal{B}(U_{n,3}) : |B \cap E_1| \leq 1\}$ and $\mathcal{B}(M_2) = \{B \in \mathcal{B}(U_{n,3}) : |B \cap E_2| \leq 2\}$. The corresponding matroids $M_1$, $M_2$ and $M_1 \cap M_2$ are of rank 3, these are given in Figure 5.

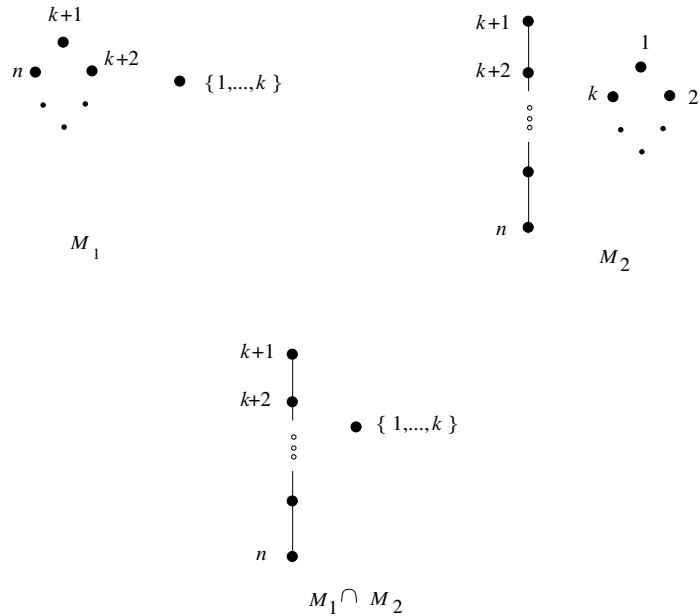

FIGURE 5. Euclidean representations of $M_1$, $M_2$ and $M_1 \cap M_2$.

2.2. **Lattice path matroids.** Let $\mathcal{A} = (A_j : j \in J)$ be a set system, that is, multiset of subsets of a finite set $E$. A *transversal* of $\mathcal{A}$ is a set $\{x_j : j \in J\}$ of $|J|$ distinct elements such that $x_j \in A_j$ for all $j$. A *partial transversal* of $\mathcal{A}$ is a transversal of a set system of the form $(A_k : k \in K)$ with $K \subseteq J$. A fundamental result due to Edmonds and Fulkerson [**?**] states that a partial transversal of a set of system $\mathcal{A} = (A_j : j \in J)$ are the independent sets of a matroid on $E$. We say that $\mathcal{A}$ is a *presentation* of the transversal matroid. The bases of a transversal matroid are the maximal partial transversals.

A *lattice path* starts at point $(0,0)$ and uses steps $(1,0)$ and $(0,1)$, called *East* and *North*. Let $P = p_1, \ldots, p_{r+m}$ and $Q = q_1, \ldots, q_{r+m}$ be two lattice paths from $(0,0)$ to $(m,r)$ with $P$ never going above $Q$. Let $\{p_{s_1}, \ldots, p_{s_r}\}$ be the set of North steps of $P$ with $s_1 < \cdots < s_r$; similarly, let $\{q_{t_1}, \ldots, q_{t_r}\}$ be the set of North steps of $Q$ with $t_1 < \cdots < t_r$. Let $M[P,Q]$ be the transversal matroid that has ground set $\{1, \ldots, m+r\}$ and presentation $(N_i : i \in \{1, \ldots, r\})$ where $N_i$ denotes the interval $[s_i, t_i]$ of integers. Transversal matroids arising as above are called *lattice path matroids*. Given a subset $X$ of $\{1, \ldots, m+r\}$, we define the lattice path $P(X) = u_1, \ldots, u_{m+r}$ where $u_i$ is a North step if $i \in X$, an East step otherwise.

In [2] was proved that a subset $B$ of $\{1, \ldots, m+r\}$ with $|B| = r$ is a base of $M[P,Q]$ if and only if the associated path $P(B)$ stays in the region bounded by $P$ and $Q$.

**Corollary 4.** *Let $M[P,Q]$ be the transversal matroid on $\{1, \ldots, m+r\}$ and presentation $(N_i : i \in \{1, \ldots, r\})$ where $N_i$ denotes the interval $[s_i, t_i]$ of integers. Suppose that there exists integer $x$ such that $s_j < x < t_j$ and $s_{j+1} < x+1 < t_{j+1}$ for some $1 \leq j \leq r-1$. Then, $P(M[P,Q])$ has a nontrivial hyperplane split.*

*Proof.* Let $E_1 = \{1, \ldots, x\}$ and $E_2 = \{x+1, \ldots, m+r\}$. Then, $M|_{E_1}$ (resp. $M|_{E_2}$) is the transversal matroid with representation $(N_i^1 : i \in \{1, \ldots, r\})$ where $N_i^1 = N_i \cap E_1$ (resp. with representation $(N_i^2 : i \in \{1, \ldots, r\})$ where $N_i^2 = N_i \cap E_2$). Let $r_1$ and $r_2$ be the ranks of $M|_{E_1}$ and $M|_{E_2}$ respectively. We have that $N_i^1 \neq \emptyset$ for all $i \leq j+1$ (since the smallest element in $N_i$ is strictly smaller than $x+1$). Therefore, $r_1 \geq j+1$. Similarly, we have that $N_i^2 \neq \emptyset$ for all $i \geq r-j+1$ (since the smallest element in $N_i$ is larger than $x+1$). Therefore, $r_2 \geq r-j+1$. So, the partition $(E_1, E_2)$ verifies property (P1) by taking integers $a_1$ and $a_2$ such that $r_1 - a_1 = j$ and $r_2 - a_2 = r-j$. Moreover, the sets $\mathcal{B}(M_1) = \{B \in \mathcal{B}(M) : |B \cap E_1| \leq r_1 - a_1\}$ and $\mathcal{B}(M_2) = \{B \in \mathcal{B}(M) : |B \cap E_2| \leq r_2 - a_2\}$ are the collections of bases of matroids $M_1$ and $M_2$ respectively. Indeed, $M_1$ is the transversal matroid with representation $(\overline{N}_i^1 : i \in \{1, \ldots, r\})$ where $\overline{N}_i^1 = N_i$ for each $i = 1, \ldots, j$ and $\overline{N}_i^1 = N_i \cap E_2$ for each $i = j+1, \ldots, r$. $M_2$ is the transversal matroid with representation $(\overline{N}_i^2 : i \in \{1, \ldots, r\})$ where $\overline{N}_i^2 = N_i \cap E_1$ for each $i = 1, \ldots, j$ and $\overline{N}_i^1 = N_i$ for each $i = j+1, \ldots, r$. Finally, $M_1 \cap M_2$ is the transversal matroid with representation $(\overline{N}_i : i \in \{1, \ldots, r\})$ where $\overline{N}_i = \overline{N}_i^1 \cap \overline{N}_i^2$ for each $i = 1, \ldots, r$. The result follows by Remark 1. $\qquad\square$

Notice that there might be several choices for the values of $x, j, a_1$ and $a_2$ (each of which arises a good partition). However, it is not clear if these partitions give different hyperplane splits. Also notice that the partition proposed in the above proof may verify neither (P2) nor (P2'). We shall see this, for instance, in the following example.

**Example 4:** Let $m = 3$ and $r = 4$. Let $P = p_1, \ldots, p_7$ be the lattice path where $p_1, p_2, p_3, p_4$ are East steps and $p_5, p_6, p_7$ are North steps . Let $Q = q_1, \ldots, q_7$ be the lattice path where $q_2, q_4, q_6$ are East steps and $q_1, q_3, q_5, q_7$ are North steps. Let $M[P, Q]$ be the transversal matroid on $\{1, \ldots, 7\}$ and presentation $(N_i : i \in \{1, \ldots, 4\})$ where $N_1 = [1, 2, 3, 4]$, $N_2 = [3, 4, 5, 6]$, $N_1 = [5, 6]$ and $N_1 = [7]$, see Figure 6.

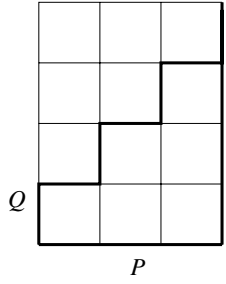

$Q$

$P$

FIGURE 6. Lattice paths $P$ and $Q$

Recall that a subset $B$ of $\{1, \ldots, 7\}$ with $|B| = 4$ is a base of $M[P, Q]$ if and only if the associated path $P(B)$ stays in the region bounded by $P$ and $Q$. So,

$$\begin{aligned}
\mathcal{B}(M) = \ & \{\{1, 3, 5, 7\}, \{1, 3, 6, 7\}, \{1, 4, 5, 7\}, \{1, 4, 6, 7\}, \{1, 5, 6, 7\}, \\
& \{2, 3, 5, 7\}, \{2, 3, 6, 7\}, \{2, 4, 5, 7\}, \{2, 4, 6, 7\}, \{2, 5, 6, 7\}, \\
& \{3, 4, 5, 7\}, \{3, 5, 6, 7\}, \{4, 5, 6, 7\}\}
\end{aligned}$$

Let us take $x = 3$ and $j = 1$, then $E_1 = \{1, 2, 3\}$ and $E_2 = \{4, 5, 6, 7\}$. So, the presentation of $M|_{E_1}$ (resp. $M|_{E_2}$) is given by $N_1^1 = [1, 2, 3]$ and $N_2^1 = [3]$ (resp. by $N_1^2 = [4, 5], N_2^2 = [4, 5, 6]$, $N_3^2 = [5, 6]$ and $N_4^2 = [7]$) and thus, $r_1 = 1$ and $r_2 = 4$. By setting $a_1 = a_2 = 1$, we obtain that $M_1$ is the transversal matroid with representation $(\overline{N}_i^1 : i \in \{1, \ldots, 4\})$ where $\overline{N}_1^1 = [1, 2, 3, 4], \overline{N}_2^1 = [4, 5, 6], \overline{N}_3^1 = [5, 6]$ and $\overline{N}_4^1 = [7]$. And, $M_2$ is the transversal matroid with representation $(\overline{N}_i^2 : i \in \{1, \ldots, 4\})$ where $\overline{N}_1^2 = [1, 2, 3], \overline{N}_2^2 = [3, 4, 5, 6], \overline{N}_3^2 = [5, 6]$ and $\overline{N}_4^2 = [7]$. Finally, we have that $M_1 \cap M_2$ is the transversal matroid with representation $(\overline{N}_i : i \in \{1, \ldots, 4\})$ where $\overline{N}_1 = [1, 2, 3], \overline{N}_2 = [4, 5], \overline{N}_3 = [5, 6]$ and $\overline{N}_4 = [7]$. These matroids are illustrated in Figure 7. We deduce that

$$\begin{aligned}
\mathcal{B}(M_1) = \ & \{\{1, 4, 5, 7\}, \{1, 4, 6, 7\}, \{1, 5, 6, 7\}, \{2, 4, 5, 7\}, \{2, 4, 6, 7\}, \\
& \{2, 5, 6, 7\}, \{3, 4, 5, 7\}, \{3, 5, 6, 7\}, \{4, 5, 6, 7\}\}, \\
\mathcal{B}(M_2) = \ & \{\{1, 3, 5, 7\}, \{1, 3, 6, 7\}, \{1, 4, 5, 7\}, \{1, 4, 6, 7\}, \{1, 5, 6, 7\}, \\
& \{2, 3, 5, 7\}, \{2, 3, 6, 7\}, \{2, 4, 5, 7\}, \{2, 4, 6, 7\}, \{2, 5, 6, 7\}, \\
& \{3, 4, 5, 7\}, \{3, 5, 6, 7\}\}
\end{aligned}$$

and

$$\begin{aligned}
\mathcal{B}(M_1 \cap \mathcal{B}(M_2)) = \ & \{\{1, 4, 5, 7\}, \{1, 4, 6, 7\}, \{1, 5, 6, 7\}, \{2, 4, 5, 7\}, \{2, 4, 6, 7\}, \\
& \{2, 5, 6, 7\}, \{3, 4, 5, 7\}, \{3, 5, 6, 7\}\}.
\end{aligned}$$

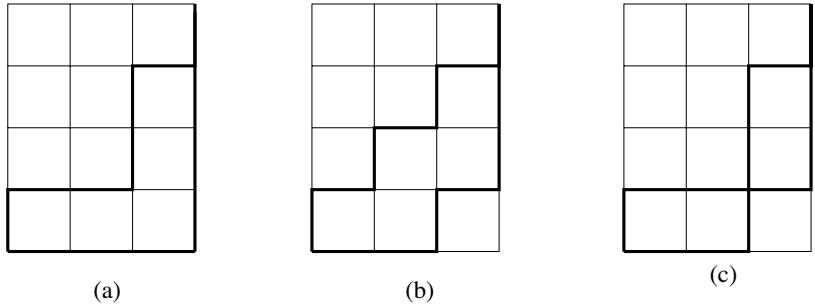

(a)                          (b)                          (c)

FIGURE 7. Transversal matroids (a) $M_1$ (b) $M_2$ and (c) $M_1 \cap M_2$

So, by Corollary 4, $M[P, Q]$ has a nontrivial hyerplane split. We notice that the above partition do not verify either (P2) (for instance, $\{3\} \in \mathcal{I}(M|_{E_1}), \{4, 5, 6\} \in \mathcal{I}(M|_{E_2})$ but $\{3, 4, 5, 6\} \notin \mathcal{I}(M[P, Q])$) or (P2') (for instance, $B = \{1, 4, 5, 7\} \in \mathcal{B}(M_1), Y = \{4, 5, 6\} \in \mathcal{I}(M|_{E_2})$ but $(B \cap E_1) \cup Y = \{1, 4, 5, 6\} \notin \mathcal{I}(M[P, Q])$).

**Example 5:** Let $m \geq r \geq 2$ be integers. Let $P_1 = p_1, \ldots, p_{m+r}$ be the lattice path where $p_1, \ldots, p_m$ are East steps and $p_{m+1}, \ldots, p_{m+r}$ are North steps. Let $Q_1 = q_1, \ldots, q_{m+r}$ be the lattice path where $q_1, \ldots q_r$ are North steps and $q_{r+1}, \ldots, q_{r+m}$ are East steps. Let $M_1[P_1, Q_1]$ be the transversal matroid on $\{1, \ldots, m+r\}$ and presentation $(N_i : i \in \{1, \ldots, r\})$ where $N_i$ denotes the interval $[s_i = i, m + i = t_i]$. So, integer $x = 2$ is such that $s_1 < x < t_1$ and $s_2 < x + 1 < t_2$ and thus, by Corollary 4, $P(M_1)$ has a nontrivial hyperplane split. Notice that any $r$-subset of $\{1, \ldots, m+r\}$ is a base of $M_1$, in other words, $M_1$ is isormorphic to $U_{m+r,r}$ (and thus, we found a particular case of Corollary 3).

We close this section by stating the following.

**Conjecture 1.** *If $P(M) = P(M_1) \cup P(M_2)$ is a nontrivial hyperplane split of $P(M)$ then the set of bases of matroids $M_1$ and $M_2$ are of the form given in Lemma 1 for a suitable partition of $E$ verifying property (P1).*

## 3. DIRECT SUM

Let $M_1 = (E_1, \mathcal{B})$ and $M_2 = (E_2, \mathcal{B})$ be matroids of rank $r_1$ and $r_2$ respectively where $E_1 \cap E_2 = \emptyset$. The *direct sum*, denoted by $M_1 \oplus M_2$, of matroids $M_1$ and $M_2$ has as ground set the disjoint union $E(M_1 \oplus M_2) = E(M_1) \cup E(M_2)$ and as set of bases $\mathcal{B}(M_1 \oplus M_2) = \{B_1 \cup B_2 | B_1 \in \mathcal{B}(M_1), B_2 \in \mathcal{B}(M_2)\}$. Further, the rank of $M_1 \oplus M_2$ is $r_1 + r_2$. Our main result in this section is the following.

**Theorem 2.** *Let $M_1 = (E_1, \mathcal{B})$ and $M_2 = (E_2, \mathcal{B})$ be matroids of rank $r_1$ and $r_2$ respectively where $E_1 \cap E_2 = \emptyset$. Then, $P(M_1 \oplus M_2)$ has a nontrivial hyperplane split if and only if either $P(M_1)$ or $P(M_2)$ has a nontrivial hyperplane split.*

We may first show the following lemma needed for the proof of Theorem 2.

**Lemma 2.** *Let $M_1 = (E_1, \mathcal{B})$ and $M_2 = (E_2, \mathcal{B})$ be matroids of rank $r_1$ and $r_2$ respectively where $E_1 \cap E_2 = \emptyset$. Then, $A \subset \mathcal{B}(M_1 \oplus M_2)$ is the collection of bases of a matroid if and only if $A = \{X \cup Y : X \in A_1, \ Y \in A_2\}$ where $A_i \subseteq \mathcal{B}(M_i)$ is the collection of bases of a matroid for each $i = 1, 2$.*

*Proof.* (*Sufficiency*) We notice that $A = \mathcal{B}(M(A_1) \oplus M(A_2))$ where $M(A_i)$ is the matroid given by the collection of bases $A_i$ on $E_i$, $i = 1, 2$. Thus, $A \subset \mathcal{B}(M_1 \oplus M_2)$ and, by definition of the direct sum, $A$ is indeed the collection of bases of a matroid.

(*Necessity*) We suppose that $A \subset \mathcal{B}(M_1 \oplus M_2)$ is the collection of bases of a matroid, say $M(A)$. We set

$$A_1 = \{B \cap E_1 : B \in A\} \text{ and } A_2 = \{B \cap E_2 : B \in A\}.$$

We first show that $A_1$ is the collection of bases of a matroid (it can also be shown that $A_2$ is the collection of bases of a matroid by using similar arguments). We thus show that the basis exchange axiom is verified. The case when $|A_1| = 1$ is clear. Let us suppose that $|A_1| \geq 2$. Let $D_1, D_2 \in A_1$. So, there exist $B_1, B_2 \in A$ such that $B_i \cap E_1 = D_i$, $i = 1, 2$. Since $M(A)$ is a matroid then if $e \in (B_1 \setminus B_2) \cap E_1 = D_1 \setminus D_2$ then there exists $f \in B_2 \setminus B_1$ such that $B_1 - e + f \in A$. Since $e \in E_1$ then $|(B_1 - e) \cap E_2| = r_2$ implying that $f \notin E_2$ (otherwise, $|(B_1 - e + f) \cap E_2| = r_2 + 1$, which is impossible). Thus, $f \in E_1$ and so $f \in D_2 \setminus D_1$. Therefore, $(B_1 - e + f) \cap E_1 = D_1 - e + f \in A_1$. We observe that the latter implies the following:

**Remark 2.** *For all $e \in D_1 \setminus D_2$ there exists $f \in D_2$ such that if $D_1 \cup Y \in A$ then $(D_1 - e + f) \cup Y \in A$ for any $Y \in \mathcal{B}(M_2)$.*

We now show that $A = \{X \cup Y : X \in A_1, \ Y \in A_2\}$. It is clear that $A \subseteq \{X \cup Y : X \in A_1, \ Y \in A_2\}$. Let $X' \in A_1$ and $Y' \in A_2$, we shall show that $X' \cup Y' \in A$ and so $\{X \cup Y : X \in A_1, \ Y \in A_2\} \subseteq A$. We first observe that for any $X' \in A_1$ and any $Y' \in A_2$ there exist $Y'' \in A_2$ and $X'' \in A_1$ such that

$$X' \cup Y'' \in A \text{ and } X'' \cup Y' \in A.$$

If $X' = X''$ then we clearly have $X' \cup Y' \in A$. Let us suppose then that $X' \neq X''$. By Remark 2, for all $e \in X'' \setminus X'$ there exists $f \in X'$ such $(X'' - e + f) \cup Y' \in A$. We can then construct a path $X'' = D_1, \ldots, D_m = X'$ connecting $X''$ to $X'$ such that $D_i \cup Y' \in A$ for each $i = 1, \ldots, m - 1$. Since, $D_1 \cup Y' = X'' \cup Y' \in A$ then we may conclude that $D_m \cup Y' = X' \cup Y' \in A$. $\square$

We may now prove Theorem 2.

*Proof of Theorem 2.* (*Necessity*) We suppose that $P(M_1 \oplus M_2) = P(A) \cup P(C)$ is a nontrivial hyperplane split for some matroids $A$ and $C$. Since $\mathcal{B}(A), \mathcal{B}(C) \subset \mathcal{B}(M_1 \oplus M_2)$ then, by Lemma 2, we have

$$\mathcal{B}(A) = \{X \cup Y : X \in A_1, Y \in A_2\} \text{ and } \mathcal{B}(C) = \{X \cup Y : X \in C_1, Y \in C_2\}$$

where $A_1, C_1 \subseteq \mathcal{B}(M_1)$ and $A_2, C_2 \subseteq \mathcal{B}(M_2)$ are the collection of bases of matroids.

We know that $\mathcal{B}(A) \cup \mathcal{B}(C) = \mathcal{B}(M_1 \oplus M_2)$. We claim that $A_2 = C_2 = \mathcal{B}(M_2)$ and $A_1, C_1 \subset \mathcal{B}(M_1)$ (or, symmetrically, $A_1 = C_1 = \mathcal{B}(M_1)$ and $A_2, C_2 \subset \mathcal{B}(M_2)$). Indeed, if $A_1 \neq \mathcal{B}(M_1)$ then there exists $X' \in \mathcal{B}(M_1)$ with $X' \notin A_1$ such that for all $Y \in \mathcal{B}(M_2)$ we have $X' \cup Y \in \mathcal{B}(A) \cup \mathcal{B}(C)$. Moreover, since $X' \cup Y \notin \mathcal{B}(A)$ then $X' \cup Y \in \mathcal{B}(C)$. Therefore, $\{X' \cup Y : Y \in \mathcal{B}(M_2)\} \subseteq \mathcal{B}(C)$ and so $C_2 = \mathcal{B}(M_2)$ implying that $C_1 \neq \mathcal{B}(M_1)$ (since $\mathcal{B}(C) \subset \mathcal{B}(M_1 \oplus M_2)$). Similarly, we may also obtain that $A_2 = \mathcal{B}(M_2)$. Now, if $A_1 = \mathcal{B}(M_1)$ then $A_2 \neq \mathcal{B}(M_2)$ (since $\mathcal{B}(A) \subset \mathcal{B}(M_1 \oplus M_2)$). We obtain, by using similar arguments as above, that $C_1 = \mathcal{B}(M_1)$ and $C_2 \neq \mathcal{B}(M_2)$.

So, suppose that $A_2 = C_2 = \mathcal{B}(M_2)$ and $A_1, C_1 \subset \mathcal{B}(M_1)$. We have that

$$\mathcal{B}(A) \cup \mathcal{B}(C) = \mathcal{B}(M_1 \oplus M_2) = \{X \cup Y : X \in A_1 \cup C_1, \ Y \in \mathcal{B}(M_2)\}$$

and thus, $A_1 \cup C_1 = \mathcal{B}(M_1)$. Also,

$$\mathcal{B}(A) \cap \mathcal{B}(C) = \{X \cup Y : X \in A_1 \cap C_1, \ Y \in \mathcal{B}(M_2)\} \subset \mathcal{B}(M_1 \oplus M_2)$$

so, since $\mathcal{B}(A) \cap \mathcal{B}(C)$ is a collection of bases of a matroid then, by Lemma 2, $A_1 \cap C_1$ is a collection of bases of a matroid. Therefore, $A_1$ and $C_1$ is a nontrivial matroid base decomposition of $M_1$.

We now show that this matroid base decomposition induce a nontrivial hyperplane split. To this end, we first show that there exists an hyperplane containing the elements corresponding to $A_1 \cap C_1$ not supporting $P(M_1)$. Let $H$ be the hyperplane corresponding to the nontrivial hyperplane split $P(M_1 \oplus M_2) = P(A) \cup P(C)$. Notice that $dim(H) \leq |E_1| + |E_2| - 2$ since $dim(P(M_1 \oplus M_2)) \leq |E_1| + |E_2|$. Without loss of generality, we may suppose that $dim(H) = |E_1| + |E_2| - 2$, otherwise we consider $P(M_1 \oplus M_2)$ embedded in $\mathbb{R}^{|E_1| + |E_2|}$ and properly extend the equation of $H$. Thus, $H$ can be written as

$$H(\mathbf{x}) = \left\{ \mathbf{x} \in \mathbb{R}^{|E_1| + |E_2|} : \langle \mathbf{e}^1, \mathbf{x} \rangle + \langle \mathbf{e}^2, \mathbf{x} \rangle = p \text{ where } \mathbf{e}^i = \sum_{j \in E_i} \alpha_j^i e_j \text{ with } \alpha_j^i \in \mathbb{R} \text{ and } p \in \mathbb{N} \right\}$$

where $e_j^i$ denotes the $j^{th}$ standard basis vector in $\mathbb{R}^{|E_i|}$ for each $i = 1, 2$.

$H$ contains the common facet of $P(A)$ and $P(C)$ (that is, the elements in $\mathcal{B}(A) \cap \mathcal{B}(C)$). Now, there exists $X' \in A_1 \cap C_1$ such that for all $Y \in \mathcal{B}(M_2)$ we have $X' \cup Y \in \mathcal{B}(A) \cap \mathcal{B}(C)$ and $\langle \mathbf{e}^1, \mathbf{x}' \rangle + \langle \mathbf{e}^2, \mathbf{y} \rangle = p$. So, for all $Y \in \mathcal{B}(M_2)$

$$\langle \mathbf{e}^2, \mathbf{y} \rangle = p - \langle \mathbf{e}^1, \mathbf{x}' \rangle = p - p' \tag{1}$$

where $p' \in \mathbb{N}$. Now, there exists $Y \in \mathcal{B}(M_2)$ such that $X \cup Y \in \mathcal{B}(A) \cap \mathcal{B}(C)$ for all $X \in A_1 \cap C_1$, so

$$\langle \mathbf{e}^1, \mathbf{x} \rangle + \langle \mathbf{e}^2, \mathbf{y} \rangle = p$$

and thus, by (1),

$$\langle \mathbf{e}^1, \mathbf{x} \rangle = p'. \tag{2}$$

Since $P(A) \cup P(C)$ is a nontrivial hyperplane split then there exists $B_1 = X_1 \cup Y_1 \in \mathcal{B}(A)$ with $X_1 \in A_1$ and $Y_2 \in \mathcal{B}(M_2)$ such that $\langle \mathbf{e}^1, \mathbf{x} \rangle + \langle \mathbf{e}^2, \mathbf{x}' \rangle > p$. Then, by (2), we have $\langle \mathbf{e}^1, \mathbf{x} \rangle + p - p' > p$ and so,

$$\langle \mathbf{e}^1, \mathbf{x} \rangle > p'. \tag{3}$$

Similarly, there exists $B_2 = X_2 \cup Y_2 \in \mathcal{B}(C)$ with $X_1 \in C_1$ and $Y_2 \in \mathcal{B}(M_2)$ such that $\langle \mathbf{e}^1, \mathbf{x} \rangle + \langle \mathbf{e}^2, \mathbf{x}' \rangle < p$. Then, by (2), we have $\langle \mathbf{e}^1, \mathbf{x} \rangle + p - p' < p$ and so,

$$\langle \mathbf{e}^1, \mathbf{x} \rangle < p'. \tag{4}$$

Therefore, by (2), the hyperplane

$$H'(\mathbf{x}) = \left\{ \mathbf{x} \in \mathbb{R}^{|E_1| + |E_2|} : \langle \mathbf{e}^1, \mathbf{x} \rangle = p' \text{ where } \mathbf{e}^1 = \sum_{j \in E_1} \alpha_j^1 e_j \text{ with } \alpha_j^1 \in \mathbb{R} \right\}$$

contains the elements corresponding to $A_1 \cap C_1$. Moreover, $H'$ does not support $P(M_1)$ by (3) and (4).

We now show that any edge of $P(M_1)$ is also and edge of either $P(A_1)$ or $P(C_1)$. We do this by contradiction, let $B_A \in A_1 \setminus C_1$ and $B_C \in C_1 \setminus A_1$ with $\Delta(B_A, B_C) = 2$. There exists $Y \in \mathcal{B}(M_2)$ such that $B_A \cup Y \in \mathcal{B}(A) \setminus \mathcal{B}(C)$ and $B_C \cup Y \in \mathcal{B}(C) \setminus \mathcal{B}(A)$. Since $\Delta(B_A \cup Y, B_C \cup Y) = \Delta(B_A, B_C) = 2$ then the edge in $P(M_1 \oplus M_2)$ joining $B_A \cup Y$ to $B_C \cup Y$ is an edge of neither $P(A)$ nor $P(C)$, which is a contradiction since $P(A) \cup P(C)$ is a nontrivial hyperplane split.

Therefore, by Proposition 1, $A_1 \cap C_1$ are the set of vertices of a common facet of $P(A_1)$ and $P(C_1)$ and, by Corollary 1, $P(M_1) = P(A_1) \cup P(C_1)$.

(*Sufficiency*) Without loss of generality, we suppose that $P(M_1) = P(N_1) \cup P(N_2)$ is a nontrivial hyperplane split for some matroids $N_i$, $i = 1, 2$. Let

$$L_1 = N_1 \oplus M_2, \ L_2 = N_2 \oplus M_2 \text{ and } L_1 \cap L_2 = (N_1 \cap N_2) \oplus M_2,$$

Since $N_1, N_2$ and $N_1 \cap N_2$ are matroids then $L_1$ and $L_2$ and $L_1 \cap L_2$ are also matroids where the collection of bases of $L_1$ and $L_2$ are

$$\mathcal{B}(L_1) = \{X \cup Y : X \in \mathcal{B}(N_1), \ Y \in \mathcal{B}(M_2)\} \text{ and } \mathcal{B}(L_2) = \{X \cup Y : X \in \mathcal{B}(N_2), \ Y \in \mathcal{B}(M_2)\}$$

Moreover, $\mathcal{B}(L_1) \cup \mathcal{B}(L_2) = \mathcal{B}(M_1 \oplus M_2)$ (since $\mathcal{B}(N_1) \cup \mathcal{B}(N_2) = \mathcal{B}(M_1)$), $\mathcal{B}(L_i) \subset \mathcal{B}(M_1 \oplus M_2)$ (since $\mathcal{B}(N_i) \subset \mathcal{B}(M_1)$ for each $i = 1, 2$) and $\mathcal{B}(L_1) \cap \mathcal{B}(L_2) \neq \emptyset$ (since $\mathcal{B}(N_1) \cap \mathcal{B}(N_2) \neq \emptyset$). Thus, the matroids given by the collection of bases $\mathcal{B}(L_1)$ and $\mathcal{B}(L_2)$ is a nontrivial matroid base decomposition of $\mathcal{B}(M_1 \oplus M_2)$.

We now show that this matroid base decomposition induces a nontrivial hyperplane split. To this end, we first show that there exists an hyperplane containing the elements in $\mathcal{B}(L_1) \cap \mathcal{B}(L_2)$ not supporting $P(M_1 \oplus M_2)$. Let $H$ be the hyperplane in $\mathbb{R}^d$ where $d$ is the dimension of $P(M_1)$ containing a common facet of $P(N_1)$ and $P(N_2)$. We suppose that $P(N_1)$ lies in the closed half-space $\overline{H}^+$ and that $P(N_2)$ lies in the other closed half-space $\overline{H}^-$ ($H$ exists since $P(N_1) \cup P(N_2)$ is a nontrivial hyperplane split). Moreover, there exist $B_1 \in \mathcal{B}(N_1)$ and $B_2 \in \mathcal{B}(N_2)$ such that $B_1$ lies in the open half-space $H^+$ while $B_2$ lies in the open half-space $H^-$. Let $H_1$ be the hyperplane defined by the same equation as $H$ in $\mathbb{R}^{d'}$ where $d'$ is the dimension of $P(M_1 \oplus M_2)$ (notice that $d' > d$). Then, $H_1$ contains the elements in $\mathcal{B}(L_1) \cap \mathcal{B}(L_2)$. Moreover, $H_1$ does not support $P(M_1 \oplus M_2)$ since $B_1 \cup Y \in \mathcal{B}(L_1)$ lies in open half-space $H_1^+$ and $B_2 \cup Y \in \mathcal{B}(L_2)$ lies in the open half-space $H_1^-$ for some $Y \in \mathcal{B}(M_2)$.

We now show that any edge of $P(M_1 \oplus M_2)$ is also and edge of either $P(L_1)$ or $P(L_2)$. We do this by contradiction, let $B_1 \in \mathcal{B}(L_1) \backslash \mathcal{B}(L_2)$ and $B_2 \in \mathcal{B}(L_2) \backslash \mathcal{B}(L_1)$ with $\Delta(B_1, B_2) = 2$. Since $B_1 = X_1 \cup X_2$ with $X_1 \in \mathcal{B}(N_1), X_1 \notin \mathcal{B}(N_2)$ and $B_2 = Y_1 \cup Y_2$ with $Y_1 \notin \mathcal{B}(N_1), Y_1 \in \mathcal{B}(N_2)$ then $X_1 \neq Y_1$. So, $X_2 = Y_2$ (since $\Delta(B_1, B_2) = 2$).

Therefore, by Corollary 1, $\mathcal{B}(L_1) \cap \mathcal{B}(L_2)$ are the set of vertices of a facet of $P(L_1)$ and $P(L_2)$ and $P(M_1 \oplus M_2) = P(L_1) \cup P(L_2)$. $\qquad \square$

**Example 6:** It is known [2] that lattice path matroids are closed under direct sums, this is illustrated in Figure 8.

Let $M[P_1, Q_1]$ and $M[P_2, Q_2]$ be two lattice path matroids. Then, by Corollary 4, $P(M[P_1, Q_1])$ (or $P(M[P_2, Q_2])$) has a nontrivial hyperplane split and so, by Theorem 2, $P(M[P_1, Q_1] \oplus M[P_2, Q_2])$ also does.

## 4. BINARY MATROIDS

Maurer [15, Theorem 2.1] gave a complete characterization of those graphs that are matroid base graphs. Let $X, Y$ be two vertices and let $\delta(X, Y)$ be the distance between $X$ and $Y$ in $G(M)$. If $\delta(X, Y) = 2$ (that is, $X$ and $Y$ are not adjacent but they are joined by a path of length two) then their *common neighbor* is defined as the set of vertices adjacent to both $X$ and $Y$. Maurer showed that each common neighbor is a square, pyramid or octahedron, see Figure 9.

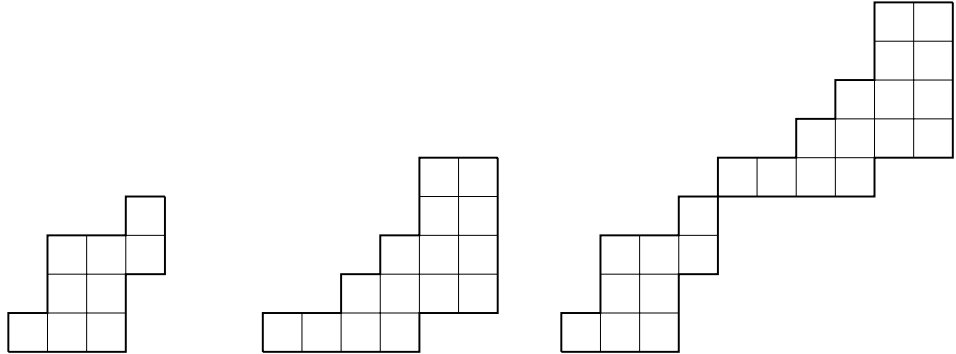

FIGURE 8. Direct sum of two lattice path matroids

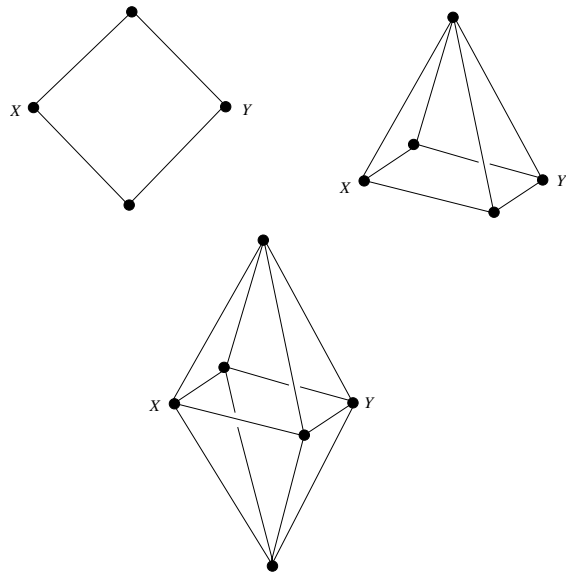

FIGURE 9. A square, a pyramid and an octahedron

In [16, Theorem 4.1] Maurer proved that a matroid $M$ is binary if and only if $G(M)$ contains no induced octahedra. Notice that $G(U_{4,2})$ is given by an octahedron (see Figure 2) and that induced octahedra in $G(M)$ correspond to *minors* of $M$ isomorphic to $U_{4,2}$.

**Corollary 5.** *Let $M = (E, \mathcal{B})$ be a binary matroid and let $X, Y \in \mathcal{B}$ with $\delta(X, Y) = 2$. Then, there exists a unique couple $U, V \in \mathcal{B}$ such that $X, U, V, Y$ form an* empty square *in $G(M)$, that is, $X, U, V, Y$ form a cycle of length four without diagonals.*

*Proof.* The existence of an empty square follows from Maurer's characterization. Indeed, for each pair of vertices $X, Y$ with $\delta(X, Y) = 2$ we have that their common neighbor is either a square or a pyramid. In both cases we can find the desired empty square. For uniqueness, suppose that there are two empty squares, say $X, U, V, Y$ and $X, U', V', Y$ with $\{U, V\} \neq \{U', V'\}$. If exactly one of $U, V$ is in $\{U', V'\}$, say $V = V'$, see Figure 10 (a) (respectively, if none of $U, V$ belongs to $\{U', V'\}$, see Figure 10 (b)) then, by the basis exchange axiom $\{X, V = V', U, U', Y\}$ must induce a pyramid, and thus neither $X, U, V, Y$ nor $X, U', V', Y$ form an empty square (respectively, $\{X, V = V', U, U', Y\}$ induce an octahedron which is not possible since $M$ is binary).

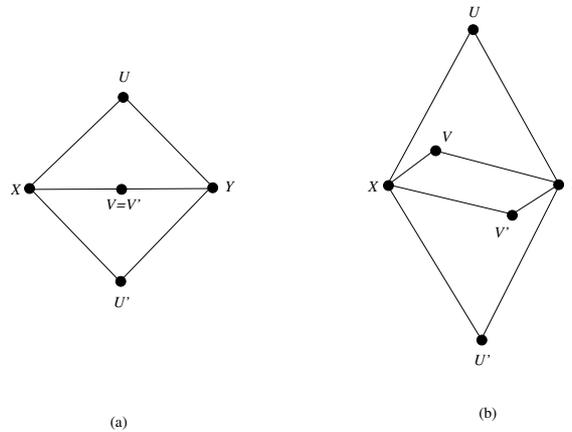

FIGURE 10. Possible empty squares

$\square$

**Lemma 3.** *Let $M = (E, \mathcal{B})$ be a binary matroid and let $\mathcal{B}_1 \subset \mathcal{B}$ such that $\mathcal{B}_1$ are the set of bases of a matroid, say $M_1$. If $X \in \mathcal{B}_1$ and the neighborhood of $X$ (that is, the set of vertices in $G(M)$ adjacent to $X$) are elements of $\mathcal{B}_1$ then $\mathcal{B}_1 = \mathcal{B}$.*

*Proof.* We first notice that $M_1$ is binary. Indeed, since $M$ is binary, we consider the columns of the matrix representation of $M$ corresponding to those elements used in $\mathcal{B}_1$. Let $X \in \mathcal{B}_1$ and let $\delta_i(X) = \{Y \in \mathcal{B}(M) | \delta(X, Y) = i\}$. We shall use induction on $i$. The cases $i = 0, 1$ are true by hypothesis. Suppose that $\bigcup_{i=1}^{k} \delta_i(X) \subset \mathcal{B}_1$ and let $Y \in \delta_{k+1}(X)$. Since $G(M)$ is connected then there exists a path joining $X$ to $Y$ and thus there is $Z \in \delta_{k-1}(X)$ such that $\delta(Z, Y) = 2$. By Corollary 5, there exists a unique couple $U, V \in \mathcal{B}$ such that $Z, U, V, Y$ form an empty square in $G(M)$. So, $\delta(X, U) = \delta(X, V) = k$ implying that $U, V \in \delta_k(X)$. Thus, since $\delta_k(X) \in \mathcal{B}_1$ and $M_1$ is binary then, by Corollary 5, there exists a unique couple $Z', Y' \in \mathcal{B}_1$ such that $Z', U, V, Y'$ form an empty square in $G(M_1)$. Since $Z', U, V, Y'$ is also an empty square in $G(M)$ then (by uniqueness) we must have that $Z' = Z$ and $Y' = Y$. Therefore, $Y \in \mathcal{B}_1$. $\qquad\square$

**Theorem 3.** *Let $M$ be a binary matroid. Then, $P(M)$ has not a nontrivial hyperplane split.*

*Proof.* By contradiction, let us suppose that $P(M) = P(M_1) \cup P(M_2)$ is a nontrivial hyperplane split for some matroids $M_1, M_2$ (and thus $P(M_1), P(M_2) \neq P(M)$). The latter induces the matroid base decomposition of $\mathcal{B}(M) = \mathcal{B}(M_1) \cup \mathcal{B}(M_2)$. By Corollary 1, $P(M_1)$ contains a vertex of $P(M)$ together with all its neighbors. But then, by Lemma 3, $M_1 = M$ which is a contradiction since the hyperplane split is nontrivial. $\qquad\square$

**Corollary 6.** *Let $M$ be a binary matroid. If $G(M)$ contains a vertex $X$ having exactly $d$ neighbours (that is, with $|\delta_1(X)| = d$) where $d = dim(P(M))$ then $P(M)$ is indecomposable.*

*Proof.* By contradiction, let $P(M) = \bigcup_{i=1}^{t} P(M_i)$, $t \geq 2$ be a decomposition. Notice that $d = dim(P(M)) = dim(P(M_i))$ for all $i$. We claim that if vertex $X$ belongs to $P(M_i)$ for some $i$ then $\delta_1(X)$ also does. Indeed, any vertex in $P(M)$ (and in a general in a $d$-polytope) must have at least $d$ neighbours since $dim(P(M_i)) = d$. So if $X$ belongs to $P(M_i)$ then all its $d$ neigbours must also belong to $P(M_i)$. Thus, since $M$ is binary and by using Lemma 3, we have that $M_i = M$ which is a contradiction since $t \geq 2$. $\qquad\square$

Maurer [17] has proved that the $d$-dimensional hypercube is the base matroid graph of a binary matroid. We have the following immediate consequence.

**Corollary 7.** *Let $P(M)$ be the base matroid polytope having as 1-skeleton the hypercube. Then, $P(M)$ is indecomposable.*

*Proof.* It follows by Corollary 6 since any vertex in the $d$-dimensional hypercube has precisely $d$ neighbours. $\qquad\square$

**Example 7:** The *Fano* matroid, denoted by $F_7$, is the rank 3 matroid on 7 elements having as set of bases all triples except $\{1, 5, 3\}, \{1, 2, 6\}, \{1, 4, 7\}, \{2, 3, 4\}, \{2, 5, 7\}, \{3, 6, 7\}$ and $\{4, 5, 6\}$, see Figure 11.

It is known that $F_7$ is binary, for instance, it can be represented by the following matrix over $\mathbb{F}^2$

$$A = \begin{pmatrix} 1 & 0 & 0 & 0 & 1 & 1 & 1 \\ 0 & 1 & 0 & 1 & 0 & 1 & 1 \\ 0 & 0 & 1 & 1 & 1 & 0 & 1 \end{pmatrix}.$$

Thus, by Theorem 3, $F_7$ has not a nontrivial hyperplane split. The following corollary is an immediate consequence of Theorem 3 since any graphic matroid is binary.

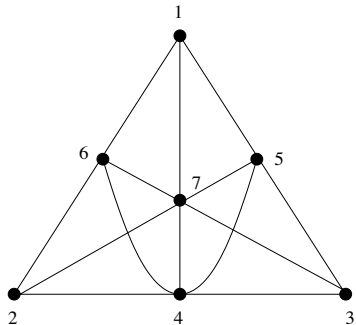

FIGURE 11. Euclidean representations of $F_7$.

**Corollary 8.** *$P(M)$ has not a nontrivial hyperplane split if $M$ is a graphic matroid.*

In particular, the matroid $M(K_4)$ has not a nontrivial hyperplane split.

Equipe Combinatoire et Optimisation, Université Pierre et Marie Curie, Paris 6, 4 Place Jussieu, 75252 Paris Cedex 05

*E-mail address*: `vchatelain,ramirez@math.jussieu.fr`